\documentclass[10pt,leqno]{amsart}
\usepackage{graphicx}
\usepackage{titlesec}

\usepackage{indentfirst,csquotes}

\usepackage{amssymb,amsthm,amsmath}
\usepackage{xcolor,paralist,hyperref,titlesec,fancyhdr,etoolbox}
\usepackage{titlesec}
\newtheorem{theorem}{Theorem} 

\newtheorem{lem}[theorem]{Lemma} 

\newtheorem{corollary}[theorem]{Corollary} 

\theoremstyle{definition} 
\theoremstyle{remark} 
\newtheorem{remark}[theorem]{Remark}

\titleformat{\section}[display]{\normalfont\huge\bfseries\centering}{\centering\chaptertitlename\thechapter}{10pt}{\Large}
\titlespacing*{\section}{2pt}{2ex}{2ex}

\hypersetup{ colorlinks=true, linkcolor=black, filecolor=black, urlcolor=black }

\usepackage{lipsum}
\begin{document}
\title{Generalizations of Some Inequalities of the Polar Derivatives of Polynomials} 
\author[Deepak Kumar, Dinesh Tripathi, Sunil Hans]{
  Deepak Kumar\textsuperscript{1}, Dinesh Tripathi\textsuperscript{2}, Sunil Hans\textsuperscript{1}
}
\date{\today}
\address{Address}
\address{Author (1): Department of Applied Mathematics, Amity University,\\ Noida-201313, India.}
\email{deepak.kumar26@s.amity.edu}
\email{Sunil.hans82@yahoo.com}
\address{Author (2): Department of Science - Mathematics, School of Sciences, Manav Rachana University, Faridabad-121004, India.}
\email{dinesh@mru.edu.in}
\let\thefootnote\relax
\begin{abstract}
The main purpose of this paper is to obtain some extensive generalizations of the polynomial inequalities in terms of the polar derivative. These generalizations depend on the consideration of the zeros inside and outside of a closed disk and the positions of extremal coefficients of the underlying polynomial.
\end{abstract} 
\maketitle
\begin{center}
\section*{INTRODUCTION}
\end{center}
The study of extremal problems in complex variables involves estimating or finding the maximum or minimum value, emphasizing inequalities related to a polynomial and its derivative under certain prescribed constraints. Usually, it applies to different classes of polynomials and distinct regions of the complex plane. These problems play a key role in geometric function theory and proving the inverse theorems in approximation theory. According to a well-known classical result due to Bernstein \cite{S}, if $P(z)$ is a polynomial of degree $n$, then
\begin{align}\label{1.1}
    \max_{|z|=1}|P'(z)| \leq n\max_{|z|=1}|P(z)|.
\end{align}
The result is the best possible, and equality in (\ref{1.1}) holds for zeros at the origin.\\
After that, in $1985,$ Frappier et al. \cite{FRC} demonstrated that if $P(z)$ is a polynomial of degree $n$, then 
\begin{align}\label{1.2}
    \max_{|z|=1}|P'(z)| \leq n\max_{1\leq l\leq 2n} |P(e^{\frac{il\pi}{n}})|.
\end{align}
It is obvious that inequality (\ref{1.2}) is a refinement of (\ref{1.1}), since the maximum value of $|P(z)|$ on $|z|=1$ may be greater than the maximum of $|P(z)|$ taken over the $(2n)^{th}$ roots of unity, as is shown by the simple example $P(z) = z^n +ia, a > 0$.\\
As a refinement of (\ref{1.2}), Aziz \cite{Az} proved that if $P(z)$ is polynomial of degree $n$, then for every real $\alpha$
\begin{align}\label{1.3}
    \max_{|z|=1} |P'(z)| \leq \frac{n}{2}(M_{\alpha}+ M_{\alpha + \pi}),
\end{align}
where
 \begin{align}\label{1.4} 
     M_\alpha = \max_{1\leq l\leq n}|P(e^{i(\alpha +2l\pi)/n})|,
     \end{align}
    and $M_{\alpha + \pi}$ is obtained by replacing $\alpha$ by $\alpha + \pi$.\\
    In the same paper, Aziz also derived that if $P(z)$ is a polynomial of degree $n$, then
    \begin{align}{\label{1.5}}
       \max_{|z|=1}|P'(z)| \geq \frac{n}{2}\biggl\{2\max_{|z|=1}|P(z)|-(M_0+M_\pi)\biggr\}. 
    \end{align}
    Inequality (\ref{1.1}) can be improved if we restrict the zeros of $P(z)$. Erd\" os conjectured and later Lax \cite{PD} proved that if $P(z)$ is a polynomial of degree $n$ having no zeros in $|z| < 1$, then
\begin{align}\label{1.6}
    \max_{|z|=1}|P'(z)| \leq \frac{n}{2}\max_{|z|=1}|P(z)|.
\end{align}
 Aziz \cite{Az} further refined inequality (\ref{1.6}) and demonstrated that if $P(z)$ is a polynomial of degree $n$ having no zeros in $|z| < 1$, then for every $\alpha$
\begin{align}\label{1.7}
    \max_{|z|=1}|P'(z)| \leq \frac{n}{2}(M^2_\alpha + M^2 _{\alpha + \pi})^\frac{1}{2},
\end{align}
where $M_\alpha$ is defined in (\ref{1.4}).\\
Nakprasit and Somsuwan \cite{NS} studied for a polynomial of degree $n$ having a zero of order $s$ at $z_0$, where $|z_0| <1$ with $\mu =1$ and proved 
\begin{theorem}{\label{thm A}}
   If $P(z)=(z-z_0)^s\sum_{v=0}^{n-s}c_vz^v$ is a polynomial of degree $n$ having no zero in $|z|<k, k\geq 1$, except a zero of order $s, 0 \leq s <n$ at $z_0$, where $|z_0|<1$, then 
   \begin{align}{\label{1.8}}
  \max_{|z|=1}|P'(z)| \leq \biggl\{\frac{s}{1-|z_0|}+\frac{A}{(1-|z_0|)^s}\biggr\}\max_{|z|=1}|P(z)|-\frac{A}{(k+|z_0|)^s}\min_{|z|=k}|P(z)|,    
   \end{align}
   where
   \[A=\frac{(1+|z_0|)^{s+1}(n-s)}{(1+k)(1-|z_0|)}.\]
\end{theorem}
Let $P(z)$ be the polynomial of degree $n$ and $\beta\in \mathbb{C}$, then 
\begin{eqnarray}\label{1.9}
	D_\beta P(z)=nP(z)+(\beta-z)P'(z)
\end{eqnarray}
is called polar derivative of polynomial $P(z)$ and the degree of $D_\beta P(z)$ is $n-1$. $D_\beta P(z)$ generalizes of ordinary derivative of $P(z)$ by dividing (\ref{1.9}) by $\beta$ and letting $\beta \to \infty$, i.e
	    \[\lim_{\beta \rightarrow\infty}\frac{D_\beta P(z)}{\beta}=P'(z).\]
Aziz \cite{A} first extended the Bernstein inequality ({\ref{1.1}}) and Erd\"{o}s-Lax result ({\ref{1.6}}) to the polar derivative. For further information about the polar derivative of polynomials, see the comprehensive books of Milovanovi\'c et al. \cite{GV}, M. Marden \cite{MM}, and Rahman and Schmeisser \cite{QG}, etc., which provide extensive resources in this direction. The latest publications, in the direction of generalization of the inequality (\ref{1.7}),  can be found in the papers \cite{SDB, JTH, DSB, WS, RC}.\\
In 2019, Mir \cite{AM} proved that if $P(z)$ is a polynomial of degree $n$ having no zeros in $|z| <1$, then for every complex number $\beta$ with $|\beta| \geq 1$,
\begin{align}{\label{1.10}}
   \max_{|z|=1}|D_\beta P(z)| \leq \frac{n}{2}\bigg[2\max_{|z|=1}|P(z)|+(|\beta|-1)\biggl\{M^2_\alpha +M^2_{\alpha+\pi}-\frac{2}{n}\left(\frac{|c_0|-|c_n|}{|c_0|+|c_n|}\right)|P(z)|\biggr\}^\frac{1}{2}\biggr],
\end{align}
where $M_\alpha$ is defined by (\ref{1.4}).\\
  As a generalization of inequality (\ref{1.10}), Hussain and Ahmad \cite{AH} derived that if  $P(z)$ is a polynomial of degree $n$ having no zeros in $|z|< k, k\leq 1$, then for every complex number $\beta$ with $|\beta| \geq 1$, 
 \begin{align}\label{1.11}
  \max_{|z|=1}|D_\beta P(z)| \leq  n\max_{|z|=1}|P(z)|+ \frac{n(|\beta|-1)}{2k^n}\Biggl\{M^2_\alpha +M^2_{\alpha +\pi}  -\frac{2k^{2n}}{n}\left(\frac{|c_0|-k^n|c_n|}{|c_0|+k^n|c_n|}\right)|P(z)|^2\Biggr\}^\frac{1}{2},
  \end{align}
     where $ M_\alpha$ is defined by $(\ref{1.4})$.\\
     As an application of inequality (\ref{1.11}), the authors \cite{AH} showed that if $P(z)=\sum_{v=0}^{n}c_v z^v, c_0 \neq 0$ is a polynomial of degree $n$ having no zeros in $|z|> k, k \geq 1$, then for every complex number $\beta$ with $|\beta | \leq 1$,
     \begin{align}{\label{1.12}}
        \max_{|z|=1}|D_\beta P(z)| \leq  n |\beta|\max_{|z|=1}|P(z)| +\frac{nk^n(1-|\beta|)}{2}\biggl\{M^2_\alpha+M^2_{\alpha+\pi} -\frac{2}{nk^n}\left(\frac{k^n|c_n|-|c_0|}{k^n|c_n|+|c_0|}\right)|P(z)|^2\biggr\}^\frac{1}{2},
     \end{align} 
      where $ M_\alpha $ is defined by $(\ref{1.4})$.\\
     Mir and Hussain \cite{MH} proved that if $P(z)=\sum_{v=0}^{n}c_v z^v, c_0 \neq 0$, with $P(z) \neq 0$ in $|z|>k, k\leq 1$, then for every complex number $\beta$ with $|\beta| \leq 1$,
     \begin{align}\label{1.13}
       \nonumber \max_{|z|=1} |D_\beta P(z)| \leq& \frac{n}{2}\Biggl[ 2|\beta|\max_{|z|=1}|P(z)|+(1-|\beta|)\Biggl\{M^2_\alpha + M^2_{\alpha + \pi}\\ &
         -\frac{2}{(1+k)}\Biggl[(1-k)+\frac{2k}{n}\Biggl(\frac{k^n|c_n|-|c_{0}|}{k^n|c_n|+|c_0|}\Biggr)\Biggr]|P(z)|^2\Biggr\}^\frac{1}{2}\Biggr],
     \end{align}
     where $ M_\alpha $ is defined by $(\ref{1.4})$.\\
     In this article, we assume that the polynomial has some zeros in the disk $|z| \geq k, k \leq 1$ or in $|z| \leq k, k \leq 1,$ or $k \geq 1$; others are outside the disk to generalize the above inequalities. Our major objective is to estimate some upper bounds for the modulus of the polar derivative of polynomials on a disk. The obtained results generalize some known estimates for the ordinary or polar derivative of polynomials as special cases.
\begin{center}
  \section*{MAIN RESULTS}	
\end{center}
  \begin{theorem}{\label{thm 2.1}}
 If $P(z)=(z-z_0)^s\sum_{v=0}^{n-s}c_vz^v$ is a polynomial of degree $n$ having no zeros in $|z| < k, k\leq 1$, except a zero of order $s$, $0\leq s < n$ at $z_0$, where $|z_0| < k, k \leq 1$, then for any complex number $\beta$ with $|\beta| \geq 1$,
 \begin{align}{\label{2.1}}
\nonumber \max_{|z|=1}|D_\beta P(z)| \leq  n \max_{|z|=1}|P(z)| & +\frac{n(|\beta|-1)}{2k^n} \Biggl\{ M^2_\alpha + M^2_{\alpha+\pi} + \frac{2s(k+|z_0|)}{n(k-|z_0|)}
 \left(\max_{|z|=1}|P(z)|\right)^2 \\&  - \frac{2k^{2n}}{n}\bigg(\frac{|c_0|-k^{n-s}|c_{n-s}|}{|c_0|+k^{n-s}|c_{n-s}|}\bigg) |P(z)|^2\Biggr\}^\frac{1}{2},
 \end{align}
 where $M_\alpha$ is defined by $(\ref{1.4})$. The result has the best possible for $ s = 0, k = 1$, and equality holds in $(\ref{2.1})$ for $P(z)= z^n+e^{i\alpha}$.
\end{theorem}
\begin{remark}
For s=0, in the Theorem \ref{thm 2.1}, we get inequality (\ref{1.11}).
If we take $z_0=0$ in Theorem \ref{thm 2.1}, we get the following result due to Mir et al. \cite{AFT}. 
\end{remark}
\begin{corollary}{\label{cor 2.1}}
If $P(z)=z^s\sum_{v=0}^{n-s}c_vz^v$ is a polynomial of degree $n$ having no zeros in $|z| < k, k\leq 1$, except a zero of order $s$, $0\leq s < n$ at the origin, then for any complex number $\beta$ with $|\beta| \geq 1$, 
\begin{align}{\label{2.2}}
   \nonumber \max_{|z|=1}|D_\beta P(z)| \leq & n \max_{|z|=1}|P(z)| +\frac{n(|\beta|-1)}{2k^n} \Biggl\{ M^2_\alpha + M^2_{\alpha+\pi} + \frac{2s}{n} \left(\max_{|z|=1}|P(z)|\right)^2 \\ & - \frac{2k^{2n}}{n}\bigg(\frac{|c_0|-k^{n-s}|c_{n-s}|}{|c_0|+k^{n-s}|c_{n-s}|}\bigg) |P(z)|^2\Biggr\}^\frac{1}{2}, 
\end{align}
 where $M_\alpha$ is defined by $(\ref{1.4})$.
  The result has the best possible for $ s = 0, k = 1$, and equality holds in $(\ref{2.2})$ for $P(z)= z^n+e^{i\alpha}$.
\end{corollary}
Dividing both the sides of inequality (\ref{2.1}) by $|\beta|$ and taking $|\beta| \to \infty$, we get the following results.
\begin{corollary}{\label{cor 2.3}}
 If $P(z)=(z-z_0)^s \sum_{v=0}^{n-s}c_vz^v$ is a polynomial of degree $n$ having no zeros in $|z| < k, k\leq 1$, except a zero of order $s$, $0\leq s < n$ at $z_0$, where $|z_0| < k, k\leq 1$, then
 \begin{align}{\label{2.3}}
    \nonumber \max_{|z|=1}|P'(z)| \leq  \frac{n}{2k^n} &\Biggl\{ M^2_\alpha + M^2_{\alpha+\pi}+\frac{2s(k+z_0|)}{n(k-|z_0|)}\left(\max_{|z|=1}|P(z)|\right)^2\\ & - \frac{2k^{2n}}{n} \bigg(\frac{|c_0|-k^{n-s}|c_{n-s}|}{|c_0|+k^{n-s}|c_{n-s}|}\bigg)|P(z)|^2\Biggr\}^\frac{1}{2},
 \end{align}
  where $M_\alpha$ is defined by $(\ref{1.4})$. The result has the best possible for $ s = 0, k = 1$, and equality holds in $(\ref{2.3})$ for $P(z)= z^n+e^{i\alpha}$.
\end{corollary} 
Next, we prove the following result as a generalization of inequality (\ref{1.13}) by applying some restriction on the zeros of $ P(z)$.  
\begin{theorem}{\label{thm 2.2}}
    If $P(z)=(z-z_0)^s\sum_{v=0}^{n-s}c_vz_v, c_0\neq 0$ is a polynomial of degree $n$ having no zeros in $|z| > k, k \leq 1$, except a zero of order $s, 0 \leq s <n $ at $z_0$, where $|z_0| > 1$, then for any complex number $\gamma$ with $|\gamma| \leq 1$, 
    \begin{align}{\label{2.4}}
       \nonumber |D_\gamma P(z)|  \leq & \frac{n}{2}\bigg[2|\gamma|\max_{|z|=1}|P(z)|+(1-|\gamma|)\biggl\{M^2_\alpha + M^2_{\alpha+\pi}-2\biggl\{1-\frac{2}{n}\biggr\{\frac{s|z_0|}{|z_0|-1}\\ &+\frac{(n-s)k}{k+1}-\frac{k}{k+1}\bigg(\frac{|a_{n-s}|k^{n-s}-|a_0|}{|a_{n-s}|k^{n-s}+|a_0|}\bigg)\biggr\}\biggr\}|P(z)|^2 \biggr\}^\frac{1}{2}\bigg],
    \end{align}
     where $M_\alpha$ is defined by $(\ref{1.4})$.  The result has the best possible for $ s = 0, k = 1$, and equality holds in $(\ref{2.4})$ for $P(z)= z^n + 1$.
\end{theorem}
\begin{remark}
    If we consider $\gamma =0$ in (\ref{2.4}), we have 
    \begin{align}{\label{2.5}}
   \nonumber |nP(z)-zP'(z)| \leq \frac{n}{2}\biggl\{M^2_\alpha & + M^2_{\alpha+\pi}-2\biggl\{1-\frac{2}{n}\biggr\{\frac{s|z_0|}{|z_0|-1}+\frac{(n-s)k}{k+1}-\\&\frac{k}{k+1}\bigg(\frac{|a_{n-s}|k^{n-s}-|a_0|}{|a_{n-s}|k^{n-s}+|a_0|}\bigg)\biggr\}\biggr\}|P(z)|^2 \biggr\}^\frac{1}{2}.    
    \end{align}
    If $\max_{|z|=1}|P(z)|=|P(e^{i\phi})|$, then we have from (\ref{2.5}) 
    \begin{align}{\label{2.6}}
    \nonumber |P'(e^{i\phi})| \geq & \frac{n}{2}\bigg[2 \max_{|z|=1}|P(z)|-\biggl\{M^2_\alpha  + M^2_{\alpha+\pi}-2\biggl\{1-\frac{2}{n}\biggr\{\frac{s|z_0|}{|z_0|-1}\\&+\frac{(n-s)k}{k+1}-\frac{k}{k+1}\bigg(\frac{|a_{n-s}|k^{n-s}-|a_0|}{|a_{n-s}|k^{n-s}+|a_0|}\bigg)\biggr\}\biggr\}|P(z)|^2 \biggr\}^\frac{1}{2}\bigg].
 \end{align}
 Since $\max_{|z|=1}|P'(z)| \geq |P'(e^{i \phi})|$, we have 
 \begin{align}{\label{2.7}}
  \nonumber \max_{|z|=1}|P'(z)| \geq & \frac{n}{2}\bigg[2 \max_{|z|=1}|P(z)|-\biggl\{M^2_\alpha  + M^2_{\alpha+\pi}-2\biggl\{1-\frac{2}{n}\biggr\{\frac{s|z_0|}{|z_0|-1}\\&+\frac{(n-s)k}{k+1}-\frac{k}{k+1}\bigg(\frac{|a_{n-s}|k^{n-s}-|a_0|}{|a_{n-s}|k^{n-s}+|a_0|}\bigg)\biggr\}\biggr\}|P(z)|^2 \biggr\}^\frac{1}{2}\bigg]. 
 \end{align}
 If we take $k=1$ in (\ref{2.7}), we get the refinement of inequality (\ref{1.5}).
Theorem \ref{thm 2.2} is reduces to inequality (\ref{1.13}) for $s=0$. 
 \end{remark}

We prove the following result as an application of Theorem \ref{thm 2.1}.
 \begin{theorem}{\label{thm 2.3}}
     If $P(z)=(z-z_0)^s\sum_{v=0}^{n-s} c_vz^v, c_0 \neq 0$ is a polynomial of degree $n$ having no zeros in $|z| > k, k \geq 1$, except a zero of order $s, 0 \leq s <n$ at $z_0$, where $|z_0| > k, k \geq 1,$ then for every complex number $\gamma$ with $|\gamma| \leq 1 $,
     \begin{align}{\label{2.8}}
       \nonumber  \max_{|z|=1}|D_\gamma P(z)| \leq & n|\gamma|\max_{|z|=1} |P(z)|+\frac{nk^n(1-|\gamma|)}{2} \biggl\{ M^2_\alpha+M^2_{\alpha+\pi} +\frac{2s(|z_0|+k)}{n(|z_0|-k)}\\& \times\left(\max_{|z|=1}|P(z)|\right)^2-\frac{2}{nk^{2n}}\left(\frac{|c_{n-s}|k^{n-s}-|c_0|}{|c_{n-s}|k^{n-s}+|c_0|}\right)|P(z)|^2 \biggr\}^\frac{1}{2}.
     \end{align}
     where $M_\alpha$ is defined by $(\ref{1.4})$.  The result has the best possible for $ s = 0, k = 1$, and equality holds in $(\ref{2.8})$ for $P(z)= z^n + 1$.
 \end{theorem} 
 \begin{remark}  
     If we take $\gamma=0$ in the above inequality, we get the following result 
     \begin{align*}
         \max_{|z|=1}|P(z)| \geq &\frac{n}{2} \bigg[2\max_{|z|=1}|P(z)|-k^n\biggl\{M^2_\alpha+M^2_{\alpha+\pi} +\frac{2s(|z_0|+k)}{n(|z_0|-k)}\\& \times\left(\max_{|z|=1}|P(z)|\right)^2-\frac{2}{nk^{2n}}\left(\frac{|c_{n-s}|k^{n-s}-|c_0|}{|c_{n-s}|k^{n-s}+|c_0|}\right)|P(z)|^2 \biggr\}^\frac{1}{2}\bigg].
     \end{align*}
      \end{remark}
      \begin{remark}
          If we take $s=0$ in Theorem \ref{thm 2.3}, we get the inequality (\ref{1.12}).
      \end{remark}
\begin{center}
\section*{LEMMAS}	
\end{center}
For the proof of our theorems, we need the following lemmas. This lemma is a special case of a result due to Govil and Rahman \cite{NK}.
\begin{lem}{\label{lem1}}
	 If $P(z)$ is a polynomial of degree $n$, and $Q(z)=z^{n}\overline{P(\frac{1}{\overline{z}})}$. Then on $|z|=1$, 
	\begin{align}{\label{3.1}}
	  |P'(z)|+|Q'(z)|\leq n \max_{|z|=1}|P(z)|.
	\end{align}
\end{lem}
The next lemma is due to Aziz \cite{Az}.
\begin{lem}{\label{lem2}}
If $P(z)$ is a polynomial of degree n. Then for $|z|=1$, and any real number $\alpha$,
\begin{align}{\label{3.2}}
    |P'(z)|^2 + |Q'(z)|^2 \leq \frac{n^2}{2}\left(M^2_\alpha +M^2_{\alpha+\pi}\right),
    \end{align}
     where $M_\alpha$ is defined by $(\ref{1.4})$.
\end{lem}
\begin{lem}{\label{lem3}}
	If $P(z)$ is a polynomial of degree $n$, then for $r \leq 1$, we have 
 \begin{align}{\label{3.3}}
    \max_{|z|=r} |P(z)| \geq r^n\max_{|z|=1}|P(z)|.
 \end{align}
\end{lem}
The above lemma is due to Varga  \cite{Varga1957}. The next lemma is due to Mir and Hussain \cite{MH}.
\begin {lem}{\label{lem4}} 
If $P(z)= \sum_{v=0}^{n} c_v z^v$ is a polynomial of degree $n$ having no zeros in $|z|< k, k\geq 1$, then for each point z on $|z|=1$ for which $P(z)\neq 0$, we have
\begin{align}{\label{3.4}}
Re\left(\frac{zP'(z)}{P(z)}\right) \leq \frac{1}{1+k}\left[n-\left(\frac{|c_0|-k^n|c_n|}{|c_0|+k^n|c_n|}\right)\right].
\end{align}
\end{lem}
\begin{lem}{\label{lem5}}
If $P(z)=(z-z_0)^s\sum_{v=0}^{n-s}c_vz^v$ is a polynomial of degree $n$ having no zeros in $|z| <k, k \geq 1$, except a zero of order $s, 0 \leq s <n$ at $z_0$, where $|z_0| <1$, then for any complex number $\beta$ with $|\beta | \geq 1$,
\begin{align}{\label{3.5}}
    \nonumber\max_{|z|=1}|D_\beta P(z)| & \leq  \frac{n}{2}\bigg[2\max_{|z|=1}|P(z)|+(|\beta|-1)\biggl\{M^2_\alpha+ M^2_{\alpha+\pi}-2\biggl\{1-\\ &\frac{2}{n}\biggl\{\frac{s}{1-|z_0|} +\frac{n-s}{1+k}-\frac{1}{1+k}\left(\frac{|a_0|-|a_n|k^{n-s}}{|a_0|+|a_n|k^{n-s}}\right)\biggr\}\biggr\}|P(z)|^2\biggr\}^\frac{1}{2}\bigg].
\end{align}
\end{lem}
The above lemma is due to I. Das et al. \cite{Das}.
\begin{center}
\section*{PROOF OF THEOREMS}	
\end{center}
\begin{proof}[\textbf{Proof of Theorem \ref{thm 2.1}}]
	Since $P(z)=(z-z_0)^s(c_0+c_1z+c_2z^2+...c_{n-s}z^{n-s})$ is a polynomial of degree $n$ having no zeros in $|z| < k,k\leq 1$, except a zero of order $s$ at $z_0$ with $|z_0|< k, k \leq 1 $.\\
    Let $F(z)= P(kz)=(kz-z_0)^sH(z)= k^s(z-\frac{z_0}{k})^s H(z)$ has no zeros in $|z| < 1$, except a zero of order $0 \leq s<n$ at $\frac{z_0}{k}, \textnormal{where } \ |\frac{z_0}{k}| <1$.
    Since, for all the $z$ on $|z|=1$, $F(z) \neq 0$, we have 
\begin{align*}
 Re\left({\frac{zF'(z)}{F(z)}}\right)= Re\left(\frac{sz}{z-\frac{z_0}{k}}\right)+Re\left(\frac{zH'(z)}{H(z)}\right).   
\end{align*}
Using Lemma \ref{lem4} for $H(z)$ of degree $(n-s)$ in the above inequality, we have 
\begin{align*}
Re\left(\frac{zF'(z)}{F(z)}\right) \leq \frac{sk}{k-|z_0|} + \frac{1}{2}\left[(n-s)-\left(\frac{|c_0|-k^{n-s}|c_{n-s}|}{|c_0|+k^{n-s}|c_{n-s}|}\right)\right].
\end{align*}
Equivalently
\begin{align}{\label{4.1}}
    Re\left(\frac{zF'(z)}{F(z)}\right) \leq \frac{n}{2}+ \frac{s(k+|z_0|)}{2(k-|z_0|)} - \frac{1}{2}\left(\frac{|c_0|-k^{n-s}|c_{n-s}|}{|c_0|+k^{n-s}|c_{n-s}|}\right).
\end{align}
If $G(z)=z^n\overline{F(\frac{1}{\overline{z}})}$, then for $|z|=1$ 
\begin{align*}
   |G'(z)|=|nF(z)-zF'(z)|.
\end{align*}
This implies for $|z|=1$
\begin{align*}
\left|\frac{zG'(z)}{F(z)}\right|^2& = \left|n-\frac{zF'(z)}{F(z)}\right|^2
 = n^2 +\left|\frac{zF'(z)}{F(z)}\right|^2-2nRe\left(\frac{zF'(z)}{F(z)}\right).
\end{align*}
\begin{align*}
    |G'(z)|^2 =n^2|F(z)|^2+|zF'(z)|^2-2nRe\left(\frac{zF'(z)}{F(z)}\right)|F(z)|^2, \ \textnormal{for} \ |z|=1.
\end{align*}
Combining inequality $(\ref{4.1})$ with above inequality for $|z|=1$, we have 
\begin{align}{\label{4.2}}
 \nonumber 2|F'(z)|^2 \leq & |F'(z)|^2 +|G'(z)|^2\\ &+2n\Bigg[\frac{s(k+|z_0|)}{2(k-|z_0|)}-\frac{1}{2}\bigg(\frac{|c_0|-k^{n-s}|c_{n-s}|}{|c_0|+k^{n-s}|c_{n-s}|}\bigg)\bigg]|F(z)|^2.   
\end{align}
Now using Lemma \ref{lem2} to (\ref{4.2}), we get
\begin{align*}
    |F'(z)|^2 \leq  \frac{n^2}{4}(V^2_\alpha + V^2_{\alpha+ \pi}) + n\Bigg[\frac{s(k+|z_0|)}{2(k-|z_0|)}-\frac{1}{2}\bigg(\frac{|c_0|-k^{n-s}|c_{n-s}|}{|c_0|+k^{n-s}|c_{n-s}|}\bigg) \bigg]|F(z)|^2,  
\end{align*}
where $V_\alpha =\displaystyle{\max_{1 \leq l \leq n}}|F(e^{i(\alpha+2l\pi)/n})|$.
\begin{align}{\label{4.3}}
   |F'(z)| \leq \frac{n}{2} \Biggl\{ V^2_\alpha + V^2_{\alpha+\pi}+\frac{2}{n} \bigg[ \frac{s(k+|z_0|)}{(k-|z_0|)}-\bigg(\frac{|c_0|-k^{n-s}|c_{n-s}|}{|c_0|+k^{n-s}|c_{n-s}|}\bigg) \bigg]|F(z)|^2\Biggr\}^\frac{1}{2}.
\end{align}
 Since $ k \leq 1$, we have by Maximum Modulus Principle
\begin{align*}
    V_\alpha=\max_{1 \leq l \leq n}|F(e^{i(\alpha+2l\pi)/n})|  = \max_{1 \leq l \leq n}|P(ke^{i(\alpha+2l\pi)/n})|  \leq \max_{1 \leq l \leq n}|P(e^{i(\alpha+2l\pi)/n})| = M_\alpha .
    \end{align*}
Therefore
\begin{align*}
    k \max_{|z|=k}|P'(z)| \leq \frac{n}{2} \Biggl\{ M^2_\alpha + M^2_{\alpha+\pi}+\frac{2}{n} \bigg[ \frac{s(k+|z_0|)}{(k-|z_0|)}-\bigg(\frac{|c_0|-k^{n-s}|c_{n-s}|}{|c_0|+k^{n-s}|c_{n-s}|}\bigg) \bigg] (\max_{|z|=k}|P(z)|)^2\Biggr\}^\frac{1}{2}.
\end{align*} 
Using Lemma \ref{lem3} to the polynomial $P'(z)$ , we get 
\begin{align*}
     k^n \max_{|z|=1}|P'(z)| \leq \frac{n}{2} \Biggl\{ M^2_\alpha + M^2_{\alpha+\pi}+\frac{2}{n} \bigg[ \frac{s(k+|z_0|)}{(k-|z_0|)}-\bigg(\frac{|c_0|-k^{n-s}|c_{n-s}|}{|c_0|+k^{n-s}|c_{n-s}|}\bigg) \bigg] (\max_{|z|=k}|P(z)|)^2\Biggr\}^\frac{1}{2}.
\end{align*}
Which gives, by the Maximum Modulus Principle for $|z|=1$ and Lemma \ref{lem3} for $P(z)$, we have 
\begin{align}{\label{4.4}}
    \nonumber \max_{|z|=1}|P'(z)| \leq  \frac{n}{2k^n} \Biggl\{& M^2_\alpha + M^2_{\alpha+\pi}+ \frac{2s(k+|z_0|)}{n(k-|z_0|)}\left(\max_{|z|=1}|P(z)|\right)^2\\&-\frac{2k^{2n}}{n}\bigg(\frac{|c_0|-k^{n-s}|c_{n-s}|}{|c_0|+k^{n-s}|c_{n-s}|}\bigg) \left(|P(z)|\right)^2\Biggr\}^\frac{1}{2}.
\end{align}
 By definition of the polar derivative of a polynomial for the complex number $ \beta$ and $|z|=1$, we have
 \begin{align}{\label{4.5}}
\nonumber |D_\beta P(z)| &= |nP(z)+(\beta-z)P'(z)| \\
\nonumber & =|nP(z)-zP'(z)+\beta P'(z)|\\
\nonumber &\leq |nP(z)-zP'(z)|+|\beta||P'(z)|\\
\nonumber &=|Q'(z)|+|P'(z)|-|P'(z)|+|\beta||P'(z)|\\
&\leq n \max_{|z|=1}|P(z)|+(|\beta|-1)|P'(z)|, \quad \text{by Lemma (\ref{lem1})}.
 \end{align}
 From (\ref{4.4}) and (\ref{4.5}), we have
 \begin{align*}
\max_{|z|=1} |D_\beta P(z)|& \leq n \max_{|z|=1}|P(z)| + \frac{n(|\beta|-1)}{2k^n}\Biggl\{ M^2_\alpha + M^2_{\alpha+\pi} +\frac{2s(k+|z_0|)}{n(k-|z_0|)}\\& \times\left(\max_{|z|=1}P(z)|\right)^2 - \frac{2k^n}{n}\bigg(\frac{|c_0|-k^{n-s}|c_{n-s}|}{|c_0|+k^{n-s}|c_{n-s}|}\bigg) \bigg]|P(z)|^2\Biggr\}^\frac{1}{2}.
 \end{align*}
 This completes the proof of Theorem \ref{thm 2.1}.
\end{proof}
\begin{proof}[\textbf{Proof of Theorem \ref{thm 2.2}}]
    Since $P(z)=(z-z_0)^s\sum_{v=0}^{n-s}c_vz^v, c_0 \neq 0$ is a polynomial of degree $n$ having no zeros in $|z| > k, k \leq 1$, except a zero of order $s$ at $z_0$, $|z_0| > 1$. Therefore, the reciprocal polynomial $Q(z)=z^n \overline{P(\frac{1}{\overline{z}})}$ has no zeros in $|z| < \frac{1}{k}, \frac{1}{k} \geq 1$, except a zero of order $s$ at $\frac{1}{\overline{z_0}}$ with  $|\frac{1}{z_0}| < 1$.\\
    Applying Lemma \ref{lem5} to the polynomial $Q(z)$, we have 
    \begin{align}{\label{4.6}}
       \nonumber |D_\beta Q(z)| \leq \frac{n}{2}\bigg[2\max_{|z|=1}|Q(z)|&+(|\beta|-1)\biggl\{V^2_\alpha+ V^2_{\alpha+ \pi}-2\biggl\{1-\frac{2}{n}\biggl\{\frac{s|z_0|}{|z_0|-1} + \\& \frac{k(n-s)}{k+1}  -\frac{k}{k+1}\bigg(\frac{|a_{n-s}|k^{n-s}-|a_0|}{|a_{n-s}|k^{n-s}+|a_0|}\bigg)\biggr\}\biggr\}|Q(z)|^2 \biggr\}^\frac{1}{2}\biggr].
    \end{align}
Since $|P(z)|=|Q(z)|$, for $|z|=1$, therefore 
\[V_\alpha= \max_{1\leq l \leq n} |Q(e^{i(\alpha+2l\pi)/n})|=\max_{1 \leq l \leq n}|P(e^{i(\alpha+2l\pi)/n})|= M_\alpha.\]
The above relation, when substituted in (\ref{4.6}), gives
\begin{align}{\label{4.7}}
\nonumber|D_\beta Q(z)| \leq \frac{n}{2}\bigg[2\max_{|z|=1}|P(z)|&+(|\beta|-1)\biggl\{M^2_\alpha+ M^2_{\alpha+ \pi}-2\biggl\{1-\frac{2}{n}\biggl\{\frac{s|z_0|}{|z_0|-1} + \\& \frac{k(n-s)}{k+1}  -\frac{k}{k+1}\bigg(\frac{|a_{n-s}|k^{n-s}-|a_0|}{|a_{n-s}|k^{n-s}+|a_0|}\bigg)\biggr\}\biggr\}|P(z)|^2 \biggr\}^\frac{1}{2}\biggr].
\end{align}
For $|z|=1$, then it is easy to verify (for example \cite{AM}, see)
\[|D_\beta Q(z)| = |\overline{\beta}||D_{1/\overline{\beta}}P(z)|.\]
Replacing $1/\overline{\beta}$ by $ \gamma$, so that $|\gamma| \leq 1$, we obtain from (\ref{4.7})
\begin{align*}
  |D_\gamma P(z)| \leq \frac{n}{2}\bigg[2|\gamma
  |\max_{|z|=1}|P(z)|&+(1-|\gamma|)\biggl\{M^2_\alpha+ M^2_{\alpha+ \pi}-2\biggl\{1-\frac{2}{n}\biggl\{\frac{s|z_0|}{|z_0|-1} + \\& \frac{k(n-s)}{k+1}  -\frac{k}{k+1}\bigg(\frac{|a_{n-s}|k^{n-s}-|a_0|}{|a_{n-s}|k^{n-s}+|a_0|}\bigg)\biggr\}\biggr\}|P(z)|^2 \biggr\}^\frac{1}{2}\biggr].  
\end{align*}
This completes the proof of Theorem \ref{thm 2.2}.
\end{proof}
\begin{proof}[\textbf{Proof of Theorem \ref{thm 2.3}}]
The proof of Theorem \ref{2.3} follows on the same lines as Theorem \ref{thm 2.2} by applying Theorem \ref{thm 2.1} instead of Lemma \ref{lem5}.  
\end{proof}
\section*{Statements and Declarations}
 \textbf{Funding:} The first author is financially supported by CSIR-UGC JRF fellowship, India (Ref. No. 221610149915 ).\\
 \textbf{Author Contribution:} All authors contributed equally to this manuscript.\\ 
\textbf{Conflict of interest:} The authors declare that they have no conflict of interest.\\
 \textbf{Data Availability:} No data was used for the research described in the article.\\
 \textbf{Ethical approval:} Not applicable.

\end{document}